\documentclass[12pt]{article}
\usepackage{graphicx}

\usepackage{amsmath,amsfonts,amssymb,amsthm,graphicx,color} \parskip1pt

\def\N{I\!\!N}

\def\s{\d s}

\def\A{\= A}
\def\a{\= a}
\def\ia{\= {a}} 

\def\u{\= u}

\def\i{\= i}
\def\N{{\it Nyancana}}

\def\suls{{\'Sulvas\u tras}}

\def\sulscu{{\'Sulvas\u tras}}
\def\sul{{\'Sulvas\u tra}}

\def\bsl{{Baudh\a yana}}
\def\asl{{\A pastamba}}
\def\ksl{{K\a ty\a yana}}
\def\msl{{M\a nava}}
\def\n{\. n}

\def\h{\d h}

\begin{document}  \date{}
\title{\ksl\ \sul\ : Some Observations}
\author{S.G. Dani}
\maketitle

\begin{abstract}

The \ksl\ \sul\ has been much less studied or discussed from a modern
perspective, even though the first English translation of two {\it
  adhy\=aya}s (chapters) from it, by Thibaut, appeared as far back as
1882. Part of the reason for this seems to be that the general
approach to the \sul\ studies has been focussed on ``the
mathematical knowledge found in them''; as the other earlier \suls,
especially \bsl\ and \asl\ substantially cover the ground in this
respect the other two \suls, \msl\ and \ksl, received much less
attention, the latter especially so. 

On the other hand the broader purpose of historical mathematical
studies  
extends far beyond cataloguing what was known in various cultures,
rather 
to understand the ethos of the respective times from a mathematical
point of view, in their own setting, in order to evolve a more complete
picture of the mathematical developments, ups as well as downs, over
history. 

Viewed from this angle, a closer look at \ksl\ \sul\ assumes
significance.  
Coming at the tail-end of the \suls\ period, after which the \suls\
tradition died down due to various historical reasons that 
are really only partly understood, makes it special in certain ways. What it
omits to mention from the body of knowledge found in the earlier \suls\ would
also be of relevance to analyse in this context, as much as what it
chooses to record. Other aspects such as the difference in language,
style, would also reflect on the context. It is the purpose here to  
explore this direction of inquiry.   
       
\end{abstract}

The performance of the {\it yajna}s in the Vedic period involved
construction of altars  
({\it vedi}) and fireplaces ({\it citi} or {\it agni}) in a variety of 
intricate shapes, such as birds, tortoise and others, of quite large
sizes; the dimensions of the {\it vedi}s often extended to over 100 feet and the {\it agni}s
could be 15 feet and 20 feet, or more, in width and length. This
warranted detailed 
description of procedures for their construction, which is the subject
of the \suls, which are parts of the {\it
  Kalpas\u tra}s associated with the {\it yajurveda}. Apart from the
direct aspect of step by step 
description in the form of manuals, the \suls\  
also include enunciation of various geometric principles involved, thereby 
setting up a body of geometric
ideas and framework. 
    
The different {\it \'s\ia kh\ia}s (branches) of the Vedic people had 
their respective versions of \sulscu, though, as may be expected, a
degree of intrinsic unity may be seen in their overall contents. 
Notwithstanding the fact that there were a large number of  {\it
  \'s\ia kh\ia}s, possibly in hundreds, only eight (or nine?) 
\suls\ with mathematical content have been known in our times.   
Baudh\a yana, \A pastamba, 
M\a nava, Maitr\a yan\i ya, Var\a ha, Saty\a \d s\a \d dha, and V\a
dul associated with the {\it K\d r\d sna yajurveda}, and 
a sole \ksl\ \sul\ associated with {\it \'sukla yajurveda} are known. 
Of these, \asl,  Var\a ha, and V\a dul\ are literally the same
\cite{Kh}. (???) \msl\ 
and  Maitr\a yan\i ya are understood to be versions of each other
(though a detailed comparison does not seem to have been made yet). 
\bsl, \asl, \msl, and \ksl\ \suls\ are independent in overall
character, even 
though, as noted above, there are many commonalities. 

The dates of the
\suls ~are uncertain; according to Kashikar, as quoted in \cite{K}, 
 the following ranges may be associated with the  composition of the
 respective \suls\ :  \bsl\ and V\a dul ~(800~-~500~BCE), \asl~and 
~\msl ~(650~-~300~BCE) and \ksl, Saty\a \d s\a \d dha and Var\a ha
~(300~BCE - 400~CE). There have also been other suggestions in this
respect (see \cite{SB} for a discussion on this), placing \bsl\ around
5th or 6th century BCE, \asl\ around 5th 
and 4th century BCE, \msl\ between them, and \ksl\ around 350 BCE. 
However all dates seem to be quite speculative, and 
there do not seem to be dependable inputs on the issue. 

The \suls
~are composed in the {\it s\u tra} (aphoristic) style,  mostly  
in prose form, though parts of some of \msl\ and \ksl\ \suls\ are
found to be in 
verse form.  The texts  have been divided by later
  commentators into convenient segments, treated as individual 
{\it s\u tras}, with numbers attached, and grouped into Chapters. As
presented in \cite{SB} 
\bsl~ has 21 Chapters adding to 285 s\u tras, \asl ~has 21 Chapters
adding to   202 s\u tras,  \msl ~has
  16 Chapters adding to 228 s\u tras, and  \ksl ~has 6 Chapters adding
  to 67 s\u tras.\footnote{An extra chapter of \ksl\ is found in the
    version given in     \cite{Kh}.}

It is an enigma in the subject that the \bsl\ \sul\ which is the oldest
happens to be the most systematic and comprehensive one in many
ways. While the 
others do have some things in addition in certain respects, there are
some crucial things that are omitted, and 
one also senses in them a general lack of harmony in the presentation. 
It is almost as if the seeds of eventual decline are embedded into the \sul\
literature, though indeed this would be a rather simplistic view to take.  

There have been many pre-modern commentaries on the \suls, and they have
proved helpful in understanding the original \suls;
(unfortunately even the dates of the commentaries can not be ascertained). 
There are commentaries of Dw\a rak\a n\a tha Yajwa and Venkate\'swara
Dik\'sit  on \bsl, of Kapardisw\a m\i, Karavinda Sw\a m\i\ and
Sundarar\a ja on \asl, and of Karka and Mah\i dhara on \ksl\ (no pre-modern
commentaries are known on \msl\ \sul.)  

The \suls\ became part of the modern global scholarship through the works of
European scholars, producing translations and edited versions of them
in European languages. George Thibaut published an English translation
of \bsl\ \sul\ 
and the commentary of Dw\a rak\a n\a tha Yajwa, in the journal {\it 
  Pandit}, published from Benaras, during 1874 - 76. Edmund Burk
brought out a translation of the  \asl\ \sul\ in German in
1901. Translation of 
\msl\ \sul\ was produced by J.M. van Gelder in 1964. There have also been
various subsequent studies of these \suls, by western as well as
Indian scholars. An English translation of
all the four \suls, with commentaries,  was brought out by S.N. Sen and
A.K. Bag~\cite{SB} in 1983.  

On the whole the \ksl\ \sul\ 
has been much less studied or discussed in modern 
writings, even though the first English translation of two {\it
  adhy\=aya}s (chapters) from it, by Thibaut, appeared as far back as
1882. Part of the reason for this seems to be that the general
approach to the \sul\ studies has been focussed on the
mathematical knowledge from the tradition as a whole; indeed, many
writers do not make adequate distinction between individual \suls\
in the overall discourse. Since the earlier \suls,
especially \bsl\ and \asl\ substantially cover the ground in this
respect the other two \suls, \msl\ and \ksl, received much less
attention, the latter especially so. On the other hand, for a fuller
understanding from a historical point of view it would be important to
study the \suls\ with attention to their individual identities,
comparisons between them etc. 

Our aim here will be to discuss the \ksl\ \sul\ in this
overall context, concentrating on its specialities in relation to the
earlier \suls, especially \bsl\ \sul. Special significance is lent to
this by the fact that \ksl\ is from  substantially 
later times, towards the fag end of the Vedic period, after which the {\it
  yajna}s lost their sheen, as a historical phenomenon, though some
feeble remnants of the idea are seen 
embodied in the later day Hindu ritual practices, including in our 
times. We shall however content ourselves with comparisons of the technical
mathematical aspects, with some modest hints of their possible
significance. The differences in various 
individual aspects would perhaps have various possible explanations,
or they could even be coincidental. However it seems worthwhile to
``identify'' them so that eventually a more comprehensive picture may
emerge,  throwing light on broader historical issues
relating to the transformation of the \sul\ literature over its
period, and in turn of the Vedic people.   

\section{The general features}

One of the most striking distinctive features of the \ksl\ \sul\ is its 
rather small size compared to the others (as the alert reader may have
noted from the figures mentioned above in this respect). What is this
economy in aid of?  Broadly speaking the difference in size is
accounted for by the fact that \ksl\ does not go into much detail of
the construction of the individual {\it citi}s, on which \bsl, for
instance, expends considerable amount of space, many chapters in terms
of the later organization of the text. \ksl\ \sul\ is largely focussed on
``theory'', though there are some parts touching upon arrangement of
the {\it vedi}s and specific features of certain {\it citi}s. Why was
such a policy adopted? It would seem that the practice of {\it yajna}s
become too 
diffuse for it to be worthwhile to go into the details, and it was
considered best to confine to discussing ``general principles''.
Notwithstanding the small size, one does not find ``stinginess'' in
the discussion of the mathematical parts. In fact in some cases there
is a lingering discourse on what may be considered quite
elementary. An example, though
rather an extreme one, is the part with the s\u tra translated
as ``Square on 
a side of 2 units is 4, on 3 units it is 9 and on 4 units it is 16.''
which is then followed by the general statement for integers and then
again separately for select fractions $\frac 12, \frac 13$ and $\frac
14$ for the length of the side.

\section{Presentation of measures}

The \bsl ~\sul~gives at the outset (in s\u tra 1.3, right after two 
s\u tras which are in the nature of an ``Abstract'' for what is to
follow) names of 
length measures with various magnitudes. There are 18 of them, the
smallest being  {\it tila} (sesame seed) is 
$\frac 1{34}$th\footnote{It was postulated by
  Thibaut 
(see \cite{T}, page 15) that the unit  owes its origin to the fact 
that they had a  formula for $\sqrt 2$ involving the fraction $\frac
1{34}$.}  of an {\it  a\n gula} (finger) and the largest one is {\it
\i \d s\a} (pole) which is 188 {\it a\n gula}s; the commonly occurring
large unit however is {\it puru\d sa} (man),
which is 120 {\it  a\n gula}s.\footnote{    
 Many of the intermediate units do not bear 
a simple fractional relation with {\it puru\d sa} however; e.g. a {\it
  b\a hu} is 36 {\it a\n gula}s, a {\it yuga} is 86 {\it a\n gula}s,
etc..}   The {\it  a\n gula} measure 
was about $\frac 34$th of an inch, or about 1.9 centimetres.
Many of the units occur infrequently, only in the description of
specific {\it vedi}s. 

\ksl\ \sul\ involves only 9 units, all from \bsl's list, except that a term
{\it vitasti} is used for a measure of 12 {\it angula}s, which is
{\it pr\a de\'sa} in \bsl, and {\it pada} is also taken to be same
measure (12 {\it angula}s). The units range only between the  {\it 
  a\n gula} to  {\it \i \d s\a}, thus excluding in particular the fine unit
measures {\it a\d nu} and {\it tila}; these are also absent in \asl\
and \msl\ \suls. 
Unlike in the \bsl\ \sul\ no systematic listing of the units  is found
in \ksl\ \sul. On the face of it 
this may suggest a greater standardisation of units over the
period, but on the other hand it could also be because \ksl\ does not
deal with many practical situations as involved in \bsl. The fact that
the small 
units do not appear in \ksl, as also in \asl\ and \msl\ \suls, 
shows that they were not part of the practical need in
the Vedic practices, and feature in \bsl\ on account of specific
theoretical preoccupations, corroborating Thibaut's hypothesis
connecting {\it tila} to a term in the expression for $\sqrt 2$. 
The theoretical inclination seems to have been lost over a period. 

\section{Cardinal directions}

Unlike the earlier \suls\ \ksl\ gives explicitly a prescription for
locating and fixing the cardinal directions, over any day. The
east-west line is obtained as the line joining the two points on a
circle drawn around the base of a pole where the shadow of the tip of
the pole falls on the circle in the course of a day. The north-south
line is then obtained  through a process of drawing a perpendicular to the east-west line; the process consisted of 
tieing the ends of a rope to two poles along the east-west line, 
stretching the rope on either side
of the line by holding it at its midpoint, and then joining the two
points marked by the midpoint on either side. 

One may wonder how the cardinal directions were
fixed in the earlier times, especially since from the beginning they
have been very 
important to the \'sulva constructions, in terms of the spiritual
motivation, and the description of the constructions has typically
been with reference 
to the east-west line. It has been suggested that it was determined by
the shadow of the pole on the equinox day, and verified by the rising
and setting points of the the star K\d rttik\a\
\cite{SB}. Interestingly, even though the method as above is
explained at the outset of \ksl\ \sul, in Chapter 7 of it there is a
verse (Chapter VII, verse 35) about determining the East as
the direction of 
rising of the stars {\it K\d rttik\a, \'Srava\d na} or {\it Pu\d sya}
or as the midway of the directions of 
rising of {\it Citr\a}  and {\it Sw\a t\i}. 
Presumably both procedures coexisted, and used
for confirmation of each other; one wonders however why then they were
not mentioned together. Notwithstanding the reasons for this and
whatever their  mutual role in practice, the 
procedure as above marks a significant advance from a broader mathematical
point of view.

\section{Construction of rectilinear figures}

Towards construction of the basic figures needed to be drawn, 
viz. rectangles, isosceles 
triangles, symmetric trapezia, with prescribed sizes, the
\suls\ principally describe the steps for drawing perpendiculars to the
line of symmetry (such a line was  a common feature of the figures
involved, it being along the east-west direction); these are
however packaged into complete procedures, as in a manual,   
for drawing the desired figures; see \cite{D2} for a discussion on
this.  

\bsl's well-known construction of the square involves the method of
drawing a perpendicular that is now a familiar compass construction in
school geometry; given a line and a point on it, at which the perpendicular is to be drawn to the line, one picks two points on the line that are equidistant from the point and located on opposite sides, and draws arcs with centres at the points with radius greater than the distance from the point - the arcs intersect in two points, one on each side of the line, and joining them provides the desired  perpendicular  to the line passing through the given point. In
this form this method is absent in \ksl, though a variation may be
said to be involved in \ksl's prescription for locating the
north-south direction, after the east-west direction is drawn, as
mentioned above. On the whole during the entire \'sulva period it was 
not common to use the compass construction as above for drawing perpendiculars, and it seems to
have disappeared by the time of \ksl. 
A method, known as \N\ (also called {\it Niranchana}) method, 
was more prevalent, and in \ksl\ it appears 
as the ``canonical'' method for drawing perpendiculars. The 
method is based on the converse of Pythagoras theorem, that in a
triangle with sides of lengths $a,b,c$
if $c^2=a^2+b^2$ then the sides with lengths $a$ and $b$ are
perpendicular to each other;\footnote{It was believed at one time that the ancient Egyptians also
  adopted such a method, but it has subsequently been discounted - see
  \cite{G}. In the case of the \suls\ however such a method is seen all over
  the place.} see \cite{D2} for details on the method and its
convenience as a tool. 
 For the $a,b, c$ as above one uses what we now call 
Pythagorean triples, the three being integers such that
$c^2=a^2+b^2$. The same of two such triples, $(3,4,5)$
and $(5,12,13)$, are involved in 
the constructions in \bsl\ and \ksl\ \suls, using {\it Nyancana}. \bsl\ includes  
a list with 5 primitive Pythagorean
triples, including the above and also  $(8,15, 17)$, $(12,35, 37)$ and 
$(7,24,25)$, though they were not adopted for use with the \N\ method
(they are noted right after the statement of the Pythagorean theorem, and presumably meant as illustrations of the theorem
- see \cite{D2}). In 
\ksl\ there is no mention of any of these other triples (or of other
new ones), though in \asl\ we find four of the above triples, excluding 
 $(7,24,25)$, used in the construction of the {\it Mah\a vedi} by the
 \N\  method. 

\section{Pythagoras theorem and its applications}

The most notable feature of the \suls ~in terms of geometric theory is
the statement of the so called Pythagoras theorem. This stands out
especially in the context of the fact that some of them, especially
\bsl, possibly predate Pythagoras. 
Actually neither the notion of a right angle nor of a right angled
triangle 
are found in the \suls, as concepts; of course right angled triangles
appeared as parts of various figures, and were implicit in the \N
~operations, but were not identified separately. Thus the statement of the
Pythagoras theorem occurs not with respect to right angled triangles,
but rather with reference to rectangles. A close translation of how it
is stated in \bsl ~would be ``the diagonal of a rectangle makes as much
(area) as (the areas) made separately by the base and the side put
together''.  The same statement also 
 appears in \ksl\ (s\u tra 2.7) where it is followed by a clause ``{\it iti
  k\s etrajn\a nam}''. The term {\it 
  k\s etra} involved in this has been translated as ``area'' by
Thibaut, but as  
``figure'' by Datta \cite{Dat}. It is argued in \cite{SB} that \suls\
use the term {\it bh\u mi} for area, so the expression as above
means ``this is the knowledge about plane figures''. Whatever be
precise nuance of the meaning, the clause is evidently intended to
emphasize the importance of the statement to the reader. In this
respect it has a pedagogical value which seems significant. 

\medskip
The Pythagoras theorem is also applied in the \suls\ for constructing 
squares with area equal to the sum of areas of
  two given  squares (including doubling of a square, called {\it
    dvikara\d ni} which is described separately), and the difference
  of areas of 
  two given squares (with unequal areas). The constructions, which are
 of course a direct application of the
Pythagoras theorem, are
  illustrated in the following Figure; see \cite{D2} for details. 

\begin{figure}[ht!]
\begin{center}
\includegraphics[width=50mm]{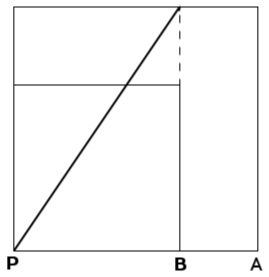}~~~~~~~
\includegraphics[width=50mm]{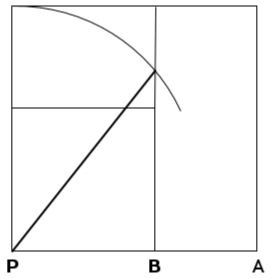}
\caption{{The thick lines give the sides of the squares with areas equal to the sum and difference,
respectively, of two given squares, with bases PA and PB as in the figure.}}
\end{center}
\end{figure}

The procedure for ``squaring'' of the difference of two squares is
also used in 
the \suls, except \msl\ \sul, for squaring a rectangle, by first
expressing it as a difference of two squares (by moving around half of the
extra part on the longer side); see \cite{D2} for details. 

The augmentation of squares was used systematically in the \suls\
for replicating
given figures in larger size, by simply enhancing the size of the reference
unit by the desired amount; e.g. to produce a replica of a figure with
area $7\frac 12$ {\it puru\d sa} to one with area $8\frac 12$ {\it
  puru\d sa}, as was required, the unit would be changed in a way that
the area will increase by a factor of $1+\frac 2{15}$; the side for
this would be obtained as a combination of the original unit with a
square 
of area $\frac 2{15}$th  of it. In fact this problem may have been the
inspiration for their discovery of the Pythagoras theorem; see
\cite{D2} for a discussion on how they may have arrived at the
theorem. 

While the conceptual framework in this respect is common to all the
\suls, \ksl\ is seen to deal with some of the
features involved more dexterously than in the earlier \suls. The use
of the method 
combining squares and rectangles seems to have become by now an art,
with the individual steps dealt with almost casually. 
In the construction of the {\it dro\d naciti} (s\u tras at the
beginning of Chapter~4) for instance, it
is quite casually prescribed to divide the square into 10 parts and to make
one of them into a square and the rest into another
square.\footnote{In \cite{SB}, both in the translation of the s\u tra
  for this (4.2 on page 123) and the commentary on it page 268, it is
  said that the square is to be divided into 100 parts; that
  interpretation is 
  however incorrect, as can be seen from \cite{Kh} and \cite{K}.}
On the problem of augmenting the unit for the purpose as described
above, \ksl\ introduces a ``rule of thumb'' method; it however does
not seem to conform to the standard stipulations accurately; apparently it was decided to pay a price in terms of accuracy in aid  of simplicity of execution. 

There is one especially notable mathematical observation found in
\ksl\ in this context. It is the procedure to produce a square with
area equal to any desired multiple of the unit. In the general \suls\
spirit this could be done by augmentation of squares of smaller
squares, starting with complete squares. \ksl\ proposes  a direct
method, that it be constructed as the altitude of an isosceles
triangle whose base is one less than the desired number (of area
multiple) and the two
equal sides add to one more than the desired number. This is an
interesting application of Pythagoras theorem and the identity 
$\frac 14 (n+1)^2-\frac 14(n-1)^2=n$. While there is indeed nothing quite like an  
abstract ``variable'' $n$ in the modern sense involved, the identity
seems to 
have been realised in terms of the desired number of unit squares to
be combined ({\it  y\a vatpram\a \d n\a ni samacaturasr\a \d ni}),
presumably by inspection of square grids.

\section{The square and the circle}
 
\ksl\ \sul\ gives the same method as \bsl\ for ``transforming'' a square into a
circle, namely for producing a circle with area equal to that of the
given square; it may be recalled that this consists of taking half the
diagonal of the square, dropping it from the centre along 
the midriff, and drawing the circle with radius which includes half
the side and a third of
the part jutting out, namely $PR$ as seen in the Figure~2 (see \cite{D2}
for details). 
 
 \begin{figure}[ht!]
\begin{center}
\includegraphics[width=65mm]{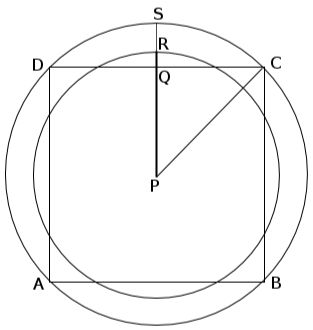}
\caption{{Circling the square: the thick line PR is given as the radius of the desired circle with area equal to that of the square ABCD.}}
\end{center}
\end{figure}

While the procedure 
is interesting (especially in its application of an intuitive
``mean-value principle'') it is not very accurate. The area of this circle
for the unit square works out to be $1.0172524...$, about $1.7 \%$ more than
that of the original square, and a computation of the value of $\pi$ 
with it comes to  $3.088 ...$). Despite this, there seems to have
been no change made from \bsl\ to \ksl. On the other hand the \msl\
\sul, though less ``sophisticated'' than either of these, and
substantially older than \ksl, seems to contain a more accurate
procedure for producing a circle with the area of a given square; see
\cite{D2}. The {Maitr\a ya\n \i ya} \sul, 
which is akin to \msl\ \sul\  gives a construction which involves
taking the radius of the desired circle to be $\frac 9{16}$ times the
side of the given square; see \cite{RCG}. Both of these involve only
about $\frac 12$ percent error. 

The converse problem of  ``squaring the circle'' viz. finding a square
with the area of a 
given circle is also considered in the \suls. Typically the treatment is not
geometric, but by assigning a numerical relation between the 
side of the desired square to the diameter of the given
circle; according to \cite{Hay}, sutra 3.2.10 from 
\msl ~is in fact a geometric construction for squaring the circle,
but we shall not go into it here; see \cite{D2}. \bsl\ \sul\
gives two formulae for squaring the circle. 
The first one gives the value for the side of the square with area
equal to that of the circle with unit radius to be $${\frac 78}+\frac
1{8\times 29}- 
\frac 1{8\times 29 \times 6}+\frac 1{8\times 29\times 6 \times 8}. \
\ \ \ \ (*)$$ 
 The area of the square with that as the side works out to be
 $3.088...$, about $98.3\%$ of the actual value. The second
 prescription consists of  taking $\frac {13}{15}$th of  
the diameter of the given circle for the side of the desired square. 
 The second is qualified by \bsl\ as
an ``incidental'' ({\it anitya}) method for squaring, signifying that
not wishing 
to use the cumbersome formula one could do with this approximate one.
Curiously, only this crude formula, giving a value that is smaller
by $4\%$ has been described in
\ksl ~\suls\ (as also \asl) for the purpose.  

\section{The square root of 2}

Three of the four \suls, 
\bsl, \asl ~and \ksl, describe a formula for $\sqrt 2$ (in
words) which corresponds to expressing the value as
$$1+\frac 13 + \frac 1{3\times 4}-\frac 1{3\times 4 \times 34}.$$
In decimal expansion the value works out to be $1.4142157\dots$, and
is noted to be accurate upto 5 decimal places.\footnote{It may be
  recalled here that the Babylonians also had a similarly close
  approximation for $\sqrt 2$, as $1.4142129\dots$, expressed in the
  sexagesimal   system that they used.} 

The s\u tra giving the formula is followed by ``{\it savi\'se\s a\h}''
in \bsl, ``{\it  sa vi\' 
  se\s a\h}'' in \asl ~and ``{\it  sa vi\'se\s a iti vi\'se\s a\h}'' in
\ksl. The word {\it vi\'se\s a\h} means ``extra''. However it does not
refer to any comparison with an accurate value of $\sqrt 2$. It is
known from the tradition of {\it \'sulvavid}s that {\it vi\'se\s a\h}
was used as a technical term for the difference between the {\it
  dvikara\d ni}, namely $\sqrt 2$, and the unit (signifying what comes
up as extra in terms of the side, while doubling a square), and in
conjunction with it  {\it savi\'se\s a\h} stood for the {\it
  dvikara\d ni} itself. Occurrence of the phrase following the s\u tra
is what tells us what the number in the s\u tra stands for (the rest of
the s\u tra only provides a number and contains no reference to what it is).  

In \cite{SB} the second part ``{\it iti vi\'se\s a\h}'' has been
translated as ``this is approximate.'' (curiously, reference the s\u
tra is missing from the commentary section in the book); see also
\cite{P}, pp. 21. However, there does not seem to be adequate 
justification for interpreting or connecting the part with an
assertion about the value being approximate. Khadilkar \cite{Kh}
translates the part, in Marathi, as ``Ha {\it dvikara\d
  ni} \d tharavi\d nyaca nira\d la prak\a ra.'', or ``This is a different
method of determining {\it dvikara\d 
  ni}''.\footnote{My translation from the Marathi version.} 
The Karka bh\a \d sya commentary
seems to confirm this; see \cite{Kh} for the commentary. The overall
context seems to favour this interpretation. From all indications the
close to accurate value 
of $\sqrt 2$ was computed, by \bsl\ or around the time, not for
practical application, for which it is not at all suitable, but to
facilitate the computation of the intricate formula for squaring the
circle,  $(*)$ as above; see \cite{Sei} and \cite{D2} for detailed
argument in this respect. However, subsequent to \bsl\ somehow no one
seems to have been interested in that formula, they being content to use
the simple proportion of 13:15 for the desired ratio. The formula for
$\sqrt 2$ then became a curiosity, bereft of its 
original significance (a nice-looking formula propounded by masters of
a bygone era). It was thus a different, unusual, formula for the
{\it dvikara\d ni}, which actually for their purpose they could simply 
measure out from the diagonal with a rope. 

There is also another interpretation possible. As noted above
 {\it vi\'se\s a\h} stands for the excess of the {\it dvikara\d
  ni} over the unit. The s\u tra thus seems to say that {\it vi\'se\s
  a\h} is such that {\it  sa vi\'se\s a} is given by the previous
expression. This is suggested by the translation in \cite{K}, where
{\it iti vi\'se\s a\h} is translated, in Hindi, as ``Yaha vi\'se\d sa
ki vy\a khy\a\ hai''. It may also be recalled here that the word
``{\it iti}'' in Sanskrit corresponds to ``in this manner'', ``thus'' or
``as you know'' (see \cite{M}) which fits well with this
interpretation.

Of course it would have been known, at least when the formula $\sqrt
2$ was first
established, that it is not exact. A formula of this kind had 
to be arrived at in some way (see \cite{D2} for a discussion on the
possibilities in this respect), and whichever way it was, it would
have been clear that an adjustment remained to be made. There is
however no indication that they considered it a significant fact worth
noting. 
\bigskip

\section{Concluding remarks} 

On the whole one notices that in the directions which were applicable
to the practical issues they met with, in terms of constructing various
rectilinear figures with conditions on the area etc., there was
progress in terms of simplifying the procedures and devising new
ones. However not much attention seems to have been paid,
collectively, to preserving 
interesting findings, even those with high aesthetic qualities, that
were not directly involved in regular practice. In some ways this could be the
result of a diffused organisation, with feeble communication and
inadequate opportunities for intellectual interaction.  

\bigskip 
\noindent{\it Acknowledgement}: The author would like to thank Kim Plofker for some useful comments on an earlier version of the manuscript and Manoj Choudhury for producing soft versions of the figures.  
 
\bigskip
\noindent 
S.G. Dani\\
UM-DAE Centre for Excellence in Basic Sciences\\
Vidyanagari Campus of University of Mumbai\\
Kalina, Mumbai 400098\\
India 

\noindent E-mail: {\tt shrigodani@cbs.ac.in}

\end{document}